\documentclass[a4paper,10pt]{article}      

\usepackage{graphicx}
\usepackage{amsmath}
\usepackage{amsfonts}
\usepackage{amssymb}
\usepackage{array,longtable}
\usepackage{amscd}
\usepackage[cp1251]{inputenc}
\usepackage[english]{babel}
\newcommand{\K}{\operatornamewithlimits{K}}
\newcommand{\TK}{\operatornamewithlimits{K^{\prime}}}
\newcommand{\low}[1]{\genfrac{}{}{0cm}{0}{}{#1}}

\newtheorem{theorem}{Theorem}[section]
\newtheorem{lemma}{Lemma}[section]
\newtheorem{corollary}{Corollary}[section]
\begin{document}
\begin{center}
{\huge {Generalized Brouncker's continued fractions and their logarithmic derivatives}} \\[0.5cm]
\end{center}

\begin{center}
    {\bf  O.Y. Kushel}
\end{center}
\begin{center}
  {\it  kushel@mail.ru}
\end{center}
\begin{center}
    Institut f\"{u}r Mathematik, MA 4-5, Technische Universit\"{a}t Berlin, \\ D-10623 Berlin, Germany
\end{center}
\begin{center}
\today
\end{center}

\begin{abstract}
In this paper, we study the continued fraction $y(s,r)$ which satisfies the equation $y(s,r)y(s+2r,r)=(s+1)(s+2r-1)$ for $r > \dfrac{1}{2}$. This continued fraction is a generalization of the Brouncker's continued fraction $b(s)$. We extend the formulas for the first and the second logarithmic derivatives of $b(s)$ to the case of $y(s,r)$. The asymptotic series for $y(s,r)$ at $\infty$ are also studied. The generalizations of some Ramanujan's formulas are presented.

Primary 11A55 Secondary 11J70, 30B70.
\end{abstract}

\section{Introduction}
\label{intro}

The Brouncker's continued fraction $b(s) = s+ \K \limits_{n=1}^{\infty}\left(\dfrac{(2n-1)^2}{2s}\right)$ still attracts the attention of researchers due to its role in the theory of orthogonal polynomials and its relations to the Gamma and Beta functions (see \cite{DUT}, \cite{KHRB}--\cite{KHR2}). Recall the following theorem of Brouncker describing the properties of $b(s)$ (see \cite{KHRB}, p. 145, Theorem 3.16).

\begin{theorem}[Brouncker]
 Let $b(s)$ be a function on $(0, + \infty)$ satisfying the functional equation $b(s)b(s+2) = (s+1)^2$ and the inequality $s < b(s)$ for $s > C$, where $C$ is some constant. Then $$b(s)
  = (s+1) \prod_{n=1}^{\infty}\dfrac{(s+4n - 3)(s+4n + 1)}{(s+ 4n - 1)^2}  =
  s +  \K_{n=1}^\infty\left(\dfrac{(2n-1)^2}{2s}\right)$$
for every positive $s$.
\end{theorem}

 Ramanujan discovered the formula, expressing the Brouncker's continued fraction in terms of the Gamma function (see \cite{KHRB}, p. 153, Theorem 3.25).

\begin{theorem}[Ramanujan] For every $s > 0$
$$b(s) = s +  \K_{n=1}^\infty\left(\dfrac{(2n-1)^2}{2s}\right) =4\left[\dfrac{\Gamma(\frac{3+s}{4})}{\Gamma(\frac{1+s}{4})}\right]^2. $$
\end{theorem}

The following extension of Brouncker's theorem (Theorem 1) was obtained by Euler
(see \cite{KHRB}, p. 180, Theorem 4.17).

\begin{theorem}[Euler] Let $y(s,r)$ be a positive continuous function satisfying the inequality $s < y(s,r)$ and the equation
$$y(s,r)y(s+2r,r) = (s+1)(s+2r-1) $$ for any $s > 0$, $r > \dfrac{1}{2}$. Then
$$y(s,r) = (s+1) \prod_{n=0}^{\infty}\dfrac{(s+2r-1+4nr)(s+4r+1+4nr)}{(s+2r+1+4nr)(s+4r-1+4nr)} = $$ $$= s+ \K_{n=1}^\infty\left(\dfrac{(2n-1)^2r^2-(r-1)^2}{2s}\right). $$
\end{theorem}

In \cite{KHRB} Ramanujan's theorem (Theorem 2) was extended to the case of the continued fraction $y(s,r)$ (see \cite{KHRB}, p. 220, ex. 4.22).

\begin{theorem}
 For every $s > 0$, $r > \frac{1}{2}$
$$y(s,r) = s+ \K_{n=1}^\infty\left(\dfrac{(2n-1)^2r^2-(r-1)^2}{2s}\right) =4r\dfrac{\Gamma(\frac{s+2r+1}{4r})\Gamma(\frac{s+4r-1}{4r})}{\Gamma(\frac{s+1}{4r})\Gamma(\frac{s+2r-1}{4r})}.$$
\end{theorem}

The following exact continued fraction representation for the first logarithmic derivative of $b(s)$
$$\dfrac{b'}{b}(s) = \dfrac{1}{s+\K\limits_{n=1}^\infty \left( \frac{n^2}{s} \right)}$$
allows one to obtain the exponential representation for $b(s)$ (see \cite{KHRB}, p. 192, Theorem 4.25).

\begin{theorem} For $s > 0$
$$ s + \K_{n=1}^\infty\left(\dfrac{(2n-1)^2}{2s}\right) = \dfrac{8\pi^2}{\Gamma^4(\frac{1}{4})}\exp \left\{ \int_0^s \dfrac{dt}{t+\K\limits_{n=1}^\infty \left( \frac{n^2}{t} \right)} \right\}.$$
\end{theorem}

In this paper, we represent the first logarithmic derivative of $y(s,r)$ in the form of the sum of two continued fractions (see Section 4, Corollary 3). For $s > |r-1|$, $r > \dfrac{1}{2}$

$$\dfrac{\partial}{\partial s}(\ln y)(s,r) = f_1(s,r) + f_2(s,r),$$

where
$$ f_1(s,r) = \dfrac{1}{2-2r+2s+2\K\limits_{n=1}^\infty \left( \frac{n^2 r^2}{1-r+s} \right)} \eqno(1)$$ and $$f_2(s,r) =\dfrac{1}{2r-2+2s+2\K\limits_{n=1}^\infty \left( \frac{n^2r^2}{r-1+s} \right)}. \eqno(2)$$

Then we extend Theorem 5 to the case of $y(s,r)$ (see Section 5, Theorem 9).

\setcounter{theorem}{8}
\begin{theorem}
 For $s > |r-1|$, $r > \frac{1}{2}$
$$ y(s,r) = s+ \K_{n=1}^\infty\left(\dfrac{(2n-1)^2r^2-(r-1)^2}{2s}\right) =$$ $$= 8\pi r2^{1 -\frac{1}{r}}\dfrac{\Gamma^2(\frac{1}{2r}) }{\Gamma^4(\frac{1}{4r})}\cot(\frac{\pi}{4r})\exp \left\{ \int_0^s (f_1(t,r) + f_2(t,r))dt \right\},$$
where $f_1(t,r)$ and $f_2(t,r)$ are given by formulas (1) and (2), respectively.\end{theorem}

There is also an exact integral representation of $\dfrac{b'}{b}(s)$ (see \cite{KHRB}, p. 191, Formula 4.71).  For $s > 0$
$$\dfrac{b'}{b}(s) = \dfrac{1}{s+\K\limits_{n=1}^\infty \left( \frac{n^2}{s} \right)} = 2\int_0^{+\infty} \dfrac{e^{-sx}dx}{\cosh x}. \eqno(3)$$

Theorem 5 together with (3) imply the following asymptotic relation, which holds for $b(s)$ as $s \rightarrow +\infty$ (see \cite{KHRB}, p. 192, Corollary 4.26).
$$b(s) = s +  \K_{n=1}^\infty\left(\dfrac{(2n-1)^2}{2s}\right) \sim
s \exp \left\{-\sum_{k=1}^{\infty}\dfrac{E_{2k}}{2ks^{2k}} \right\},$$
where $E_{2k}$ are the Euler's numbers. Here the asymptotic power series \linebreak $\exp \left\{-\sum\limits_{k=1}^{\infty}\dfrac{E_{2k}}{2ks^{2k}} \right\}$ arises from replacing $x$ in the formal power series \linebreak $\exp(x) = \sum\limits_{n=0}^\infty \dfrac{x^n}{n!}$ by $-\sum\limits_{k=1}^{\infty}\dfrac{E_{2k}}{2ks^{2k}}$ and combining coefficients afterwards (on the possibility of such a substitution see \cite{SIR}, p. 15, Theorem 124, see also \cite{DE}, p. 15).

We obtain exact integral representations for both the continued fractions (1) and (2) (see Section 3, Lemma 4). For $s> |r-1| $, $r>\frac{1}{2}$
 $$\dfrac{1}{2-2r+2s + 2\K \limits_{n=1}^\infty \left(\frac{n^2r^2}{1-r+s}\right)} = \frac{1}{2r}\int_0^{+ \infty}\dfrac{e^{-x\frac{1-r+s}{r}}dx}{\cosh x}.$$
$$\dfrac{1}{2r-2+2s + \K\limits_{n=1}^\infty \left(\frac{n^2r^2}{r-1+s}\right)} = \frac{1}{2r}\int_0^{+ \infty }\dfrac{e^{-x\frac{r-1+s}{r}}dx}{\cosh x}.$$

These two formulas together with Theorem 9 allows us to obtain the asymptotic expansion for $y(s,r)$ at infinity (see Section 6, Theorem 10), using Euler's methods.

$$y(s,r) = s+ \K\limits_{n=1}^\infty\left(\dfrac{(2n-1)^2r^2-(r-1)^2}{2s}\right) \sim $$ $$\sim  s \exp \left\{- \sum_{n=1}^{\infty}\dfrac{\sum\limits_{k = 0}^{n}\binom{2n}{2k} (r-1)^{2k} r^{2(n-k)} E_{2(n-k)} }{2ns^{2n}}\right\}. $$

Let us introduce the notation:

$$s^2-1 + \dfrac{4 \times 1^2}{1} \low{+} \dfrac{4 \times 1^2}{s^2-1} \low{+} \dfrac{4 \times 2^2}{1} \low{+} \dfrac{4 \times 2^2}{s^2-1} \low{+} \low{\cdots}
= s^2 -1 + \TK_{n=1}^{\infty}\left(\frac{4n^2}{1} \low{+} \frac{4n^2}{s^2-1}\right). $$

Ramanujan stated the following formula for the second logarithmic derivative of $b(s)$, which was proved later by Perron (see \cite{KHRB}, p. 231, Formula (5.6), see also \cite{PER}).

\setcounter{theorem}{5}
\begin{theorem} [Ramanujan's formula]
 For $s>1$
$$(\ln b)''(s) = - \int_0^{\infty} \dfrac{xe^{-sx}}{\cosh x}dx = - \dfrac{1}{s^2 -1 + \TK\limits_{n=1}^{\infty}\left(\frac{4n^2}{1} \low{+} \frac{4n^2}{s^2-1}\right)}. \eqno(4)$$
\end{theorem}

We obtain the corresponding formula for the second logarithmic derivative of $y(s,r)$.

\setcounter{theorem}{11}
\begin{theorem} For $s> \max(1, \ 2r-1)$, $r> \dfrac{1}{2}$
$$  \dfrac{\partial^2}{\partial s^2}(\ln y)(s,r) =  - \dfrac{1}{2r^2}\int_0^{\infty} \dfrac{x(e^{-\frac{1-r+s}{r}x}+ e^{-\frac{r-1+s}{r}x})}{\cosh x}dx = - h_1(s,r) - h_2(s,r), $$
where
$$h_1(s,r) = \dfrac{1}{2(1-2r+s)(1+s) + 2\TK\limits_{n=1}^{\infty}\left(\frac{4n^2r^2}{1} \low{+} \frac{4n^2r^2}{(1-2r+s)(1+s)}\right)},$$
$$ h_2(s,r) = \dfrac{1}{2(2r-1+s)(s-1) + 2\TK\limits_{n=1}^{\infty}\left(\frac{4n^2r^2}{1} \low{+} \frac{4n^2r^2}{(2r-1+s)(s-1)}\right)}.$$
\end{theorem}

\section{Functional equations for logarithmic derivatives of $y(s,r)$}
\label{sec:1} \setcounter{theorem}{6}

Let us recall the following statement, which will be used later (see \cite{KHRB}, p. 152, Lemma 3.23).

\begin{lemma} Let $g(s)$ be a monotonic function on $(0, \infty)$, vanishing at infinity, and $a > 0$ be a positive number. Then the functional equation
$f(s) + f(s+a) = g(s)$ has a unique solution, vanishing at infinity, given by the formula $$f(s) = \sum_{n=0}^\infty(-1)^n g(s+na).$$ \end{lemma}

Let us prove two following statements for the first and the second logarithmic derivatives of $y(s,r)$.

\begin{lemma} The functional equation
$$f(s,r)+f(s+2r,r)=\dfrac{1}{s+1}+\dfrac{1}{s+2r-1} = \dfrac{2(s+r)}{(s+1)(s+2r-1)} \eqno(5)$$
has a unique solution, satisfying $ \lim\limits_{s \to \infty}f(s,r)=0$, which is $$f(s,r) = \dfrac{\partial}{\partial s}(\ln y)(s,r).$$
\end{lemma}

{\bf Proof}. The equality $y(s,r)y(s+2r,r) = (s+1)(s+2r-1)$ implies
$$\ln(y(s,r)y(s+2r,r)) = \ln((s+1)(s+2r-1)); $$
$$\ln y(s,r)+ \ln y(s+2r,r) = \ln(s+1) + \ln(s+2r-1)). $$
Differentiating by $s$, we obtain
$$ \dfrac{\frac{\partial}{\partial s}y(s,r)}{y(s,r)} + \dfrac{\frac{\partial}{\partial s}y(s+2r,r)}{y(s+2r,r)} = \dfrac{1}{s+1}+\dfrac{1}{s+2r-1}. \eqno(6)$$
The function $f(s) = \dfrac{\frac{\partial}{\partial s}y}{y}(s,r)$ satisfy the conditions of Lemma 1 with $a = 2r$ and $g(s) = \dfrac{2(s+r)}{(s+1)(s+2r-1)}$. Applying Lemma 1, we complete the proof.\hfill $\square$

Let us examine two equations:

$$f_1(s,r) + f_1(s+2r,r) = \dfrac{1}{s+1} \eqno(7)$$
and
$$f_2(s,r) + f_2(s+2r,r) = \dfrac{1}{s+2r-1}. \eqno(8) $$

Both of them satisfy the conditions of Lemma 1 with $a = 2r$, $g(s) = \dfrac{1}{s+1} $ and $g(s) = \dfrac{1}{s+2r-1}$, respectively. So, applying Lemma 1, we obtain, that the solution $f_1(s,r)$ of equation (7) which satisfies $\lim\limits_{s \rightarrow \infty}f_1(s,r) = 0$ is unique. The solution $f_2(s,r)$ of equation (8) which satisfies $\lim\limits_{s \rightarrow \infty}f_2(s,r) = 0$ is also unique. Since their sum $f_1(s,r) + f_2(s,r)$ satisfies equation (5), we have from Lemma 2, that $$\dfrac{\partial}{\partial s}(\ln y)(s,r)
= f_1(s,r) + f_2(s,r),$$ where $f_1(s,r)$ and $f_2(s,r)$ are the solutions of (7) and (8), respectively, vanishing as $s \rightarrow +\infty$.

\begin{lemma} The functional equation
$$f(s,r)+f(s+2r,r)=- \dfrac{1}{(s+1)^2}-\dfrac{1}{(s+2r-1)^2} \eqno(9)$$
has a unique solution, satisfying $ \lim\limits_{s \to \infty}f(s,r)=0$, which is $$f(s,r) = \dfrac{\partial^2}{\partial^2 s}(\ln y)(s,r).$$
\end{lemma}

{\bf Proof}. Differentiate equation (6) once again by $s$:
$$ \dfrac{\partial}{\partial s}\left(\dfrac{\frac{\partial}{\partial s}y(s,r)}{y(s,r)}\right) + \dfrac{\partial}{\partial s}\left(\dfrac{\frac{\partial}{\partial s}y(s+2r,r)}{y(s+2r,r)}\right) = - \dfrac{1}{(s+1)^2}-\dfrac{1}{(s+2r-1)^2} .$$
Applying Lemma 1 with $a = 2r$ and $g(s) = -\dfrac{1}{(s+1)^2}-\dfrac{1}{(s+2r-1)^2}$, we complete the proof.\hfill $\square$

Repeating the above reasoning, we obtain that $$- \dfrac{\partial ^2}{\partial^2 s}(\ln y)(s,r)
= h_1(s,r) + h_2(s,r), \eqno(10)$$ where $h_1(s)$ is the unique solution of the equation $$h_1(s,r)+h_1(s+2r,r)=\dfrac{1}{(s+1)^2}, \eqno(11)$$ satisfying $\lim\limits_{s \to \infty}h_1(s,r)=0$ and $h_2(s)$ is the unique solution of the equation $$h_2(s,r)+h_2(s+2r,r)=\dfrac{1}{(s+2r-1)^2}, \eqno(12)$$ satisfying $\lim\limits_{s \to \infty}h_2(s,r)=0$.

\section{Exact integral representation for certain type continued fractions}
\label{sec:3}

To begin, we formulate the following result by Euler (see \cite{EU} and \cite{KHRB}, p. 191, Theorem 4.24).

\begin{theorem} For $s >0$
$$\dfrac{1}{s + \K\limits_{n=1}^\infty \left(\frac{n^2}{s}\right)} = 2\int_0^1\dfrac{x^s dx}{1+x^2}. \eqno(13)$$
\end{theorem}

\begin{corollary} For $s > 0$
$$\dfrac{1}{s+\K\limits_{n=1}^\infty \left( \frac{n^2}{s} \right)} = \int_0^{+\infty} \dfrac{e^{-sx}dx}{\cosh x}. \eqno(14)$$
\end{corollary}

Let us formulate and prove the following lemma.

\begin{lemma} Let $\varphi(s,r)$ be an arbitrary real-valued function of $s $ and $r$. Then for $r>0$, $\varphi(s,r) > 0$\footnote{Actually we have the condition $\dfrac{\varphi(s,r)}{r} > 0$ but since in the conditions of Theorem 3 $r > \dfrac{1}{2}$ we restrict ourselves to the case $r>0$.}
$$\dfrac{1}{2\varphi(s,r) + 2\K\limits_{n=1}^\infty \left(\frac{n^2r^2}{\varphi(s,r)}\right)} = \frac{1}{r}\int_0^1\dfrac{x^{\frac{\varphi(r,s)}{r}}dx}{1+x^2}. \eqno(15)$$
\end{lemma}

{\bf Proof}.
Examine equality (13). Using the substitution $s:= \dfrac{\varphi(s,r)}{r}$, where $\varphi(s,r)$ is an arbitrary real-valued function of $s$ and $r$, we obtain the equality, which is correct for all $s, \ r$ satisfying $\varphi(s,r) > 0$, $r > 0$.

$$\dfrac{1}{\frac{\varphi(s,r)}{r} + \K\limits_{n=1}^\infty \left(\frac{n^2}{\frac{\varphi(s,r)}{r}}\right)} = 2\int_0^1\dfrac{x^{\frac{\varphi(s,r)}{r}} dx}{1+x^2}. $$

$$\dfrac{1}{2\varphi(s,r) + 2r\K\limits_{n=1}^\infty \left(\frac{n^2}{\frac{\varphi(s,r)}{r}}\right)} = \dfrac{1}{r}\int_0^1\dfrac{x^{\frac{\varphi(s,r)}{r}} dx}{1+x^2}. $$

Let us apply the equivalence transform with the parameters $r_0 = 1$, $r_n = r$, $n = 1,2, \ldots$ to the continued fraction on the left-hand side. This results the formula:

$$\dfrac{1}{2\varphi(s,r) + 2\K\limits_{n=1}^\infty \left(\frac{n^2r^2}{\varphi(s,r)}\right)} = \dfrac{1}{r}\int_0^1\dfrac{x^{\frac{\varphi(s,r)}{r}} dx}{1+x^2}. $$ \hfill $\square$

\begin{corollary}
For $r>0$, $\varphi(s,r) > 0$
$$\dfrac{1}{\varphi(s,r) + \K\limits_{n=1}^\infty \left(\frac{n^2r^2}{\varphi(s,r)}\right)} = \frac{1}{r}\int_0^{+\infty} \dfrac{e^{-\frac{\varphi(r,s)}{r}x}dx}{\cosh x}.$$
\end{corollary}

{\bf Proof}.
It is enough for the proof to use the substitution $x:= e^{-x}$.\hfill $\square$

{\bf Example.}
Let $\varphi(s,r) = s+ \sin r$, $s = 1$, $r = \dfrac{\pi}{2}$. Then
$$\dfrac{1}{2 + \K\limits_{n=1}^\infty \left(\frac{n^2\frac{\pi^2}{4}}{2}\right)} = \frac{2}{\pi}\int_0^{+\infty} \dfrac{e^{-\frac{4}{\pi}x}dx}{\cosh x}.$$
Using the equivalence transformation with the parameters $r_0 = 1$, $r_n = 2$, $n = 1, \ 2, \ \ldots$, we get
$$ \dfrac{1}{4 + \K\limits_{n=1}^\infty \left(\frac{n^2\pi^2}{4}\right)} = \frac{1}{\pi}\int_0^{+\infty} \dfrac{e^{-\frac{4}{\pi}x}dx}{\cosh x}.$$

\section{Functional equations for certain type continued fractions}
\label{sec:4}

Now let us formulate and prove the following theorem.

\begin{theorem}
Let $\varphi(s,r) = s + \psi(r)$, where $\psi(r)$ is an arbitrary real-valued function of $r$. Then for $r>0$, $s> - \psi(r)$ the continued fraction of the form $$f(s,r) = \dfrac{1}{2\varphi(s,r) + 2\K\limits_{n=1}^\infty \left(\frac{n^2r^2}{\varphi(s,r)}\right)}$$ is the unique solution of the functional equation $$f(s,r) + f(s+2r,r) = \dfrac{1}{\varphi(s,r) + r}, \eqno(16)$$
satisfying $\lim\limits_{s \rightarrow \infty} f(s,r) = 0$.
\end{theorem}

{\bf Proof}.
Examine the series expansion for the right-hand side of equation (15):
$$\frac{1}{r}\int_0^1\dfrac{x^{\frac{\varphi(s,r)}{r}}dx}{1+x^2} = \dfrac{1}{r}\sum_{n=0}^{\infty}(-1)^n \int_0^1 x^{2n + \frac{\varphi(s,r)}{r}}dx =
\dfrac{1}{r}\sum_{n=0}^{\infty}\dfrac{(-1)^n}{\frac{\varphi(s,r)}{r}+2n+1} =$$
$$=\sum_{n=0}^{\infty}\dfrac{(-1)^n}{\varphi(s,r)+2rn+r} .$$

It follows from Lemma 1, that the unique solution of equation (16) satisfying $\lim\limits_{s \rightarrow \infty} f(s,r) = 0$ is given by the formula:
$$f(s,r) = \sum_{n=0}^{\infty}\dfrac{(-1)^n}{\varphi(s+2rn,r)+r}. $$
Since $\varphi(s,r) = s + \psi(r)$, we have $$f(s,r) = \sum_{n=0}^{\infty}\dfrac{(-1)^n}{\varphi(s+2rn,r)+r} = \sum_{n=0}^{\infty}\dfrac{(-1)^n}{\varphi(s,r)+2rn +r} = $$
 $$= \frac{1}{r}\int_0^1\dfrac{x^{\frac{\varphi(s,r)}{r}}dx}{1+x^2} = \dfrac{1}{2\varphi(s,r) + 2\K\limits_{n=1}^\infty \left(\frac{n^2r^2}{\varphi(s,r)}\right)},$$ by Lemma 4. \hfill $\square$

\begin{corollary} For $s> |r-1|$, $r > 0$ functional equation (7) has a unique solution satisfying $\lim\limits_{s \rightarrow 0}f_1(s,r) = 0$ which is
$$f_1(s,r) = \dfrac{1}{2-2r+2s + 2\K\limits_{n=1}^\infty \left(\frac{n^2r^2}{1-r+s}\right)}. \eqno(17)$$
Functional equation (8) also has a unique solution satisfying $\lim\limits_{s \rightarrow 0}f_2(s,r) = 0$ which is
$$f_2(s,r) = \dfrac{1}{2r-2+2s + 2\K\limits_{n=1}^\infty \left(\frac{n^2r^2}{r-1+s}\right)}. \eqno(18)$$
\end{corollary}
{\bf Proof}.
From the equality $\dfrac{1}{s+1} = \dfrac{1}{\varphi_1(s,r)+r}$ we obtain $\varphi_1(s,r) = s+1-r$. Since $\varphi_1(s,r)$ satisfies the conditions of Theorem 8, we obtain, that the continued fraction $f_1(s,r)$ is the solution of (7).
By analogy, from the equality $\dfrac{1}{s+2r-1} = \dfrac{1}{\varphi_2(s,r)+r}$ we obtain $\varphi_2(s,r) = s+r-1$, which also satisfies the conditions of Theorem 8. Applying Theorem 8 again, we obtain that the continued fraction $f_2(s,r)$ is the solution of (8).\hfill $\square$

\section{The exponential formula for generalized Brouncker's continued fraction}
\label{sec:5}

\begin{theorem}
 For $s > |r-1|$, $r > \dfrac{1}{2}$
$$ y(s,r) = s+  \K_{n=1}^\infty\left(\dfrac{(2n-1)^2 r^2-(r-1)^2}{2s}\right) =$$ $$= 8\pi r2^{1 -\frac{1}{r}}\dfrac{\Gamma^2(\frac{1}{2r}) }{\Gamma^4(\frac{1}{4r})}\cot(\frac{\pi}{4r})\exp \left\{ \int_0^s f_1(t,r)dt \right\}\exp \left\{ \int_0^s f_2(t,r)dt \right\},$$ where the continued fractions $f_1(s,r)$ and $f_2(s,r)$ are defined by equations (17) and (18), respectively. \end{theorem}

{\bf Proof}.
According to Corollary 3, the continued fractions $f_1(s,r)$ and $f_2(s,r)$ satisfy equations (7) and (8), respectively. Hence applying Lemma 2 we obtain, that
$$ \dfrac{\partial}{\partial s}(\ln y)(s,r) = \dfrac{1}{2-2r+2s + 2\K\limits_{n=1}^\infty \left(\frac{n^2r^2}{1-r+s}\right)} + \dfrac{1}{2r-2+2s + 2\K\limits_{n=1}^\infty \left(\frac{n^2r^2}{r-1+s}\right)}.$$
Integrating the obtained differential equation, we get $$\ln y(s,r) = \int_{0}^{s}(f_1(t,r)+f_2(t,r))dt + C(r)$$
$$y(s,r) = C(r)\exp\left\{\int_{0}^{s}(f_1(t,r)+f_2(t,r))dt \right\},$$
where $C(r)$ is a function of $r$.

 It is easy to see, that $C(r) = y(0,r)$. Let us calculate $y(0,r)$, using Theorem 4.
At first let us recall some well-known formulas for the Gamma function. These are the diplication formula
$$\Gamma(z)\Gamma(z+ \frac{1}{2})= 2^{1-2z}\sqrt{\pi}\Gamma(2z),$$
and the Euler's reflection formula
$$\Gamma(1-z)\Gamma(z) = \dfrac{\pi}{\sin(\pi z)}.$$
Since by definition $\Gamma(z+1) = z \Gamma(z)$, we have $\Gamma(1-z)= -z\Gamma(-z)$ and rewrite the Euler's reflection formula in the following form
$$\Gamma(z)\Gamma(-z)=  -\dfrac{\pi}{z\sin(\pi z)}. \eqno(19)$$

Using simple calculations we obtain, that
$$y(0,r)= 4r\dfrac{\Gamma(\frac{2r+1}{4r})\Gamma(\frac{4r-1}{4r})}{\Gamma(\frac{1}{4r})\Gamma(\frac{2r-1}{4r})} = 4r\dfrac{\Gamma(\frac{1}{4r} + \frac{1}{2})\Gamma(1 -\frac{1}{4r})}{\Gamma(\frac{1}{4r})\Gamma(\frac{1}{2} - \frac{1}{4r})} = \ldots$$
Since $$\Gamma(\frac{1}{2} + \frac{1}{4r})\Gamma(\frac{1}{4r}) = 2^{1-\frac{1}{2r}}\sqrt{\pi}\Gamma(\frac{1}{2r})$$ and $$\Gamma( \frac{1}{2}- \frac{1}{4r})\Gamma(- \frac{1}{4r}) = 2^{1+\frac{1}{2r}}\sqrt{\pi}\Gamma(- \frac{1}{2r}),$$ we have
$$\ldots = 4r\dfrac{2^{1-\frac{1}{2r}}\sqrt{\pi}\Gamma(\frac{1}{2r})\Gamma(1 -\frac{1}{4r})\Gamma(- \frac{1}{4r})}{\Gamma^2(\frac{1}{4r})2^{1+\frac{1}{2r}}\sqrt{\pi}\Gamma(- \frac{1}{2r})}= 4r\dfrac{\Gamma(\frac{1}{2r})\Gamma(1 -\frac{1}{4r})\Gamma(\frac{1}{4r})\Gamma(- \frac{1}{4r})}{\Gamma^3(\frac{1}{4r})2^{\frac{1}{r}}\Gamma(- \frac{1}{2r})} = $$ $$=
4r\dfrac{\Gamma(\frac{1}{2r})\pi\Gamma(- \frac{1}{4r})}{\Gamma^3(\frac{1}{4r})2^{\frac{1}{r}}\Gamma(- \frac{1}{2r})\sin( \frac{\pi}{4r})} = \ldots $$ Using (19) we obtain
$$\ldots = 4r\dfrac{\Gamma(\frac{1}{2r})\pi\Gamma(- \frac{1}{4r})\Gamma(\frac{1}{4r})}{\Gamma^4(\frac{1}{4r})2^{\frac{1}{r}}\Gamma(- \frac{1}{2r})\sin( \frac{\pi}{4r})} = - 4r\dfrac{\Gamma(\frac{1}{2r})\pi^2 }{\Gamma^4(\frac{1}{4r})2^{\frac{1}{r}}\Gamma(- \frac{1}{2r})\sin( \frac{\pi}{4r})\frac{1}{4r}\sin(\frac{\pi}{4r})} =$$
$$= - 16r^2\dfrac{\Gamma^2(\frac{1}{2r})\pi^2 }{\Gamma^4(\frac{1}{4r})2^{\frac{1}{r}}\Gamma(- \frac{1}{2r})\Gamma( \frac{1}{2r})\sin( \frac{\pi}{4r})\sin(\frac{\pi}{4r})} =  16r^2\dfrac{\Gamma^2(\frac{1}{2r})\pi^2 \frac{1}{2r}\sin(\frac{\pi}{2r})}{\Gamma^4(\frac{1}{4r})2^{\frac{1}{r}}\pi\sin( \frac{\pi}{4r})\sin(\frac{\pi}{4r})} =$$ $$= 8\pi r\dfrac{\Gamma^2(\frac{1}{2r}) 2\sin(\frac{\pi}{4r})\cos(\frac{\pi}{4r})}{\Gamma^4(\frac{1}{4r})2^{\frac{1}{r}}\sin( \frac{\pi}{4r})\sin(\frac{\pi}{4r})} = 8\pi r\dfrac{\Gamma^2(\frac{1}{2r}) }{\Gamma^4(\frac{1}{4r})2^{\frac{1}{r}-1}}\cot(\frac{\pi}{4r}).$$\hfill $\square$

\begin{corollary} For $s > 0$
$$ s + \K_{n=1}^\infty\left(\dfrac{(2n-1)^2}{2s}\right) = \dfrac{8\pi^2}{\Gamma^4(\frac{1}{4})}\exp \left\{ \int_0^s \dfrac{dt}{t+\K\limits_{n=1}^\infty \left( \frac{n^2}{t} \right)} \right\}.$$
\end{corollary}

{\bf Proof}. Just put $r=1$ and observe that
$$ \dfrac{1}{2s+2\K\limits_{n=1}^\infty \left( \frac{n^2}{s} \right)} + \dfrac{1}{2s+2\K\limits_{n=1}^\infty \left( \frac{n^2}{s} \right)} = \dfrac{2}{2s+2\K\limits_{n=1}^\infty \left( \frac{n^2}{s} \right)}= \dfrac{1}{s+\K\limits_{n=1}^\infty \left( \frac{n^2}{s} \right)}.$$\hfill $\square$

{\bf Example.}
Putting $r = 2$ into the statement of Theorem 9 and calculating $\cot(\dfrac{\pi}{8}) = \dfrac{\sin(\dfrac{\pi}{4})}{1-\cos(\dfrac{\pi}{8})} = \sqrt{2} + 1$, we obtain for $s> 1$
$$ s+  \K_{n=1}^\infty\left(\dfrac{4(2n-1)^2-1}{2s}\right) = 16\pi(2+\sqrt{2}) \dfrac{\Gamma^2(\frac{1}{4}) }{\Gamma^4(\frac{1}{8})} \times $$ $$\times \exp \left\{ \int_0^s \dfrac{dt}{2t-2+2\K\limits_{n=1}^\infty \left( \frac{4n^2}{t-1} \right)} \right\}\exp \left\{ \int_0^s \dfrac{dt}{2t+2+2\K\limits_{n=1}^\infty \left( \frac{4n^2}{t+1} \right)} \right\}.$$

\section{Generalized Brouncker's continued fraction and its asymptotic series}
\label{sec:6}
Let us recall the following lemma (see \cite{ANR}, p. 614, also \cite{KHRB}, p. 150, Lemma 3.21).

\begin{lemma}[Watson]
Let $f$ be a function on $(0, + \infty)$, such that $|f(t)| < M$ for $t > \epsilon$ and $f(t) = \sum \limits_{k=0}^{\infty}c_k t^k$, $0 < t < 2\epsilon$. Then
$$\int_0^{+ \infty}f(t)e^{-st}dt \sim \sum_{k=0}^{\infty}\dfrac{k!c_k}{s^{k+1}}, \qquad s \rightarrow + \infty $$ is the asymptotic expansion for the Laplace transform of $f$.
\end{lemma}

Let us write the asymptotic expansions for both continued fractions (17) and (18).
Applying Corollary 2, we get the following formulas for $s> |r-1|$, $r>0$:
 $$\dfrac{1}{2-2r+2s + 2\K\limits_{n=1}^\infty \left(\frac{n^2r^2}{1-r+s}\right)} = \frac{1}{2r}\int_0^{+ \infty}\dfrac{e^{-x\frac{1-r+s}{r}}dx}{\cosh x}; \eqno(20)$$
$$\dfrac{1}{2r-2+2s + 2\K\limits_{n=1}^\infty \left(\frac{n^2r^2}{r-1+s}\right)} = \frac{1}{2r}\int_0^{+ \infty }\dfrac{e^{-x\frac{r-1+s}{r}}dx}{\cosh x}. \eqno(21)$$

Examine equation (20). Write the right-hand side of equation (20) in the following form:
$$\frac{1}{2r}\int_0^{+ \infty}\dfrac{e^{-x\frac{1-r+s}{r}}dx}{\cosh x} =\frac{1}{2r}\int_0^{+ \infty} e^{\frac{r-1}{r}x} \dfrac{1}{\cosh x}e^{-\frac{s}{r}x}dx. \eqno(22)$$
Repeating the reasoning from \cite{KHRB}, p. 92, we obtain:
$$\dfrac{1}{\cosh x} = \sum_{n = 0}^{\infty}\dfrac{E_n}{n!}x^n, $$
where $E_n$ are the Euler's numbers;
$$e^{\frac{r-1}{r}x} = \sum_{n = 0}^{\infty}\dfrac{(r-1)^n}{r^n n!}x^n. $$

Using the rules of series multiplication, we get:
$$ \dfrac{ e^{\frac{r-1}{r}x}}{\cosh x} = (\sum_{n = 0}^{\infty}\dfrac{(r-1)^n}{r^n n!}x^n)(\sum_{n = 0}^{\infty}\dfrac{E_n}{n!}x^n) = \sum_{n = 0}^{\infty}\left(\sum_{k = 0}^{n}\dfrac{(r-1)^k}{r^k k!}\dfrac{E_{n-k}}{(n-k)!} \right)x^n.$$

Applying Watson's lemma 5 to (22) with $f(x) = \dfrac{e^{\frac{r-1}{r}x}}{\cosh x}$, we obtain:

$$\frac{1}{2r}\int_0^{+ \infty}\dfrac{e^{-x\frac{1-r+s}{r}}dx}{\cosh x} \sim \frac{1}{2r}\sum_{n=0}^{\infty}n!\left(\sum_{k = 0}^{n}\dfrac{(r-1)^k}{r^k k!}\dfrac{E_{n-k}}{(n-k)!} \right)\dfrac{r^{n+1}}{s^{n+1}}$$ as $s \rightarrow \infty$.

Since $\dfrac{n!}{k!(n-k)!} = \binom{n}k$, we have

$$\dfrac{1}{2-2r+2s + 2\K\limits_{n=1}^\infty \left(\frac{n^2r^2}{1-r+s}\right)} \sim   \frac{1}{2}\sum_{n=0}^{\infty}\dfrac{\sum\limits_{k = 0}^{n}\binom{n}k (r-1)^k r^{n-k} E_{n-k} }{s^{n+1}} \eqno(23)$$ as $s \rightarrow \infty$.

Analogically, we obtain the following asymptotic expansion for (21):

$$\dfrac{1}{2r - 2 +2s + 2\K\limits_{n=1}^\infty \left(\frac{n^2r^2}{r-1+s}\right)} \sim \frac{1}{2}\sum_{n=0}^{\infty}\dfrac{\sum\limits_{k = 0}^{n}\binom{n}k (1 - r)^k r^{n-k} E_{n-k} }{s^{n+1}} \eqno(24)$$ as $s \rightarrow \infty$.

\begin{theorem}
The following asymptotic relation holds as $s \rightarrow + \infty$:
$$s+ \K_{n=1}^\infty\left(\dfrac{(2n-1)^2r^2-(r-1)^2}{2s}\right) \sim $$ $$ \sim s \exp \left\{- \sum_{n=1}^{\infty}\dfrac{\sum\limits_{k = 0}^{n}\binom{2n}{2k} (r-1)^{2k} r^{2(n-k)} E_{2(n-k)} }{2ns^{2n}}\right\}. \eqno(25)$$
\end{theorem}

{\bf Proof}.
By Theorem 3, the left-hand side of (25) is divisible by $(s+1)$. Theorem 9 implies, that the continued fraction $y(s,r)$ can be written as
$$y(s,r) = (s+1)y(0,r) \exp\left\{\int_0^{+ \infty}\gamma_1(t,r)dt \right\}\exp\left\{-\int_s^{+ \infty}\gamma_1(t,r)dt \right\}\times$$ $$\times\exp\left\{\int_0^{+ \infty}\gamma_2(t,r)dt \right\}\exp\left\{-\int_s^{+ \infty}\gamma_2(t,r)dt \right\}, $$
where $$\gamma_1(t,r) = \dfrac{1}{2-2r+2t + 2\K\limits_{n=1}^\infty \left(\frac{n^2r^2}{1-r+t}\right)} - \dfrac{1}{2(1+t)},$$
$$\gamma_2(t,r) = \dfrac{1}{2r-2+2t + 2\K\limits_{n=1}^\infty \left(\frac{n^2r^2}{r-1+t}\right)} - \dfrac{1}{2(1+t)}.$$
Using asymptotic expansions (23), (24) and the expansion
$$\dfrac{1}{(1+t)} = \dfrac{1}{t(1+ \frac{1}{t})} \sim \dfrac{1}{t} \sum_{n=0}^{\infty}(-1)^n\frac{1}{t^n} =  \sum_{n=0}^{\infty}(-1)^n\frac{1}{t^{n+1}} \qquad t \rightarrow + \infty$$
 we obtain
$$\gamma_1(t,r) \sim \frac{1}{2}\sum_{n=0}^{\infty}\dfrac{\sum\limits_{k = 0}^{n}\binom{n}k (r-1)^k r^{n-k} E_{n-k} }{t^{n+1}} - \frac{1}{2}\sum_{n=0}^{\infty}(-1)^n\frac{1}{t^{n+1}} = \ldots.$$  Since the numerator of the null's term in the first sum is equal to $E_0 = 1$,
$$ \ldots =  \frac{1}{2}\sum_{n=1}^{\infty}\dfrac{\sum\limits_{k = 0}^{n}\binom{n}k (r-1)^k r^{n-k} E_{n-k} - (-1)^n}{t^{n+1}}. \qquad t \rightarrow + \infty$$
Analogically,
$$\gamma_2(t,r) \sim  \frac{1}{2}\sum_{n=1}^{\infty}\dfrac{\sum\limits_{k = 0}^{n}\binom{n}k (1-r)^k r^{n-k} E_{n-k} - (-1)^n}{t^{n+1}}. \qquad t \rightarrow + \infty$$

Integrating this over $(s, \ + \infty)$, we obtain
$$\int_s^{+ \infty} \gamma_1(t,r)dt \sim \frac{1}{2}\sum_{n=1}^{\infty}\dfrac{\sum\limits_{k = 0}^{n}\binom{n}k (r-1)^k r^{n-k} E_{n-k} - (-1)^n}{ns^{n}}. \qquad t \rightarrow + \infty$$
$$\int_s^{+ \infty} \gamma_2(t,r)dt \sim  \frac{1}{2}\sum_{n=1}^{\infty}\dfrac{\sum\limits_{k = 0}^{n}\binom{n}k (1-r)^k r^{n-k} E_{n-k} - (-1)^n}{ns^{n}}. \qquad t \rightarrow + \infty$$

Since $y(s,r) \sim s$ as $s \rightarrow + \infty$, we conclude that
$$ y(0,r) \exp\left\{\int_0^{+ \infty}\gamma_1(t,r)dt \right\}\exp\left\{\int_0^{+ \infty}\gamma_2(t,r)dt \right\} = 1$$
and
$$y(s,r) \sim (s+1)\exp\left\{- \frac{1}{2}\sum_{n=1}^{\infty}\dfrac{\sum\limits_{k = 0}^{n}\binom{n}k (r-1)^k r^{n-k} E_{n-k} - (-1)^n}{ns^{n}} \right\} \times $$ $$ \times \exp\left\{-\frac{1}{2}\sum_{n=1}^{\infty}\dfrac{\sum\limits_{k = 0}^{n}\binom{n}k (1-r)^k r^{n-k} E_{n-k} - (-1)^n}{ns^{n}} \right\}.$$

Using the equality $\sum\limits_{n=1}^{\infty} \dfrac{(-1)^n}{ns^n} = -\ln\left(\dfrac{s+1}{s}\right),$ as $s>1$, we obtain that
$$y(s,r) \sim s\exp\left\{- \frac{1}{2}\sum_{n=1}^{\infty}\dfrac{\sum\limits_{k = 0}^{n}\binom{n}k ((r-1)^k + (1-r)^k) r^{n-k} E_{n-k}}{ns^{n}} \right\} = $$
$$= \ [\mbox{ since} (1-r)^k = (-1)^k(r-1)^k] \ =$$
$$= s\exp\left\{- \sum_{n=1}^{\infty}\dfrac{\sum\limits_{k = 0}^{[\frac{n}{2}]}\binom{n}{2k} (r-1)^{2k} r^{n-2k} E_{n-2k}}{ns^{n}} \right\}.$$
The proof is completed by observing that all the Euler's numbers with odd parameters $E_1$, $E_3$, $E_5$, $\ldots$ are equal to zero.\hfill $\square$

{\bf Example.}
Putting $r = 2$ we obtain
$$s+ \K_{n=1}^\infty\left(\dfrac{4(2n-1)^2-1}{2s}\right) \sim s \exp \left\{- \sum_{n=1}^{\infty}\dfrac{\sum\limits_{k = 0}^{n}\binom{2n}{2k}2^{2(n-k)} E_{2(n-k)} }{2ns^{2n}}\right\},$$
as $s \rightarrow + \infty$.
Computations with the first few Euler's numbers $E_0 = 1$, $E_1 = 0$, $E_2 = -1$, $E_3 = 0$, $E_4 = 5$, $E_5 = 0$, $E_6 = -61$ shows that $\sum\limits_{k=0}^1\binom{2}{2k}2^{2(1-k)} E_{2(1-k)} = -3$ for $n = 1$, $\sum\limits_{k=0}^2\binom{4}{2k}2^{2(2-k)} E_{2(2-k)} = 57$ for $n = 2$ and $\sum\limits_{k=0}^3\binom{6}{2k}2^{2(3-k)} E_{2(3-k)} = -2763$ for $n=3$.

So we have $$s+ \K_{n=1}^\infty\left(\dfrac{4(2n-1)^2-1}{2s}\right) \sim s \exp \left\{\dfrac{3}{2s^2} -\dfrac{57}{4s^4} + \dfrac{2763}{6s^6} + O\left(\dfrac{1}{s^8}\right)\right\}.$$

Writing the first terms of the expansion of $e^x$
$$e^x \sim 1+x+\dfrac{x^2}{2}+\dfrac{x^3}{6} + O(x^4)$$ and substituting $x = \dfrac{3}{2s^2} -\dfrac{57}{4s^4} + \dfrac{2763}{6s^6} + O\left(\dfrac{1}{s^8}\right)$ we obtain the first terms of the expansion:
$$s+ \K_{n=1}^\infty\left(\dfrac{4(2n-1)^2-1}{2s}\right) \sim \left(1+ \dfrac{3}{2s^2} - \dfrac{105}{8s^4} + \dfrac{7035}{16s^6} + O\left(\dfrac{1}{s^8}\right) \right) = $$ $$= s+ \dfrac{3}{2s} - \dfrac{105}{8s^3} + \dfrac{7035}{16s^5} + O\left(\dfrac{1}{s^7}\right).$$

\section{Ramanujan's formula and its generalization}
\label{sec:7}

Our generalization of Ramanujan's formula (4) requires some preliminary results. The first of them is the following theorem.
\begin{theorem}
Let $\varphi(s,r)$ be an arbitrary real-valued function of $s$ and $r$. Then for $r>0$, $\varphi(s,r) > r$
$$\dfrac{1}{\varphi(s,r) - r^2 + \TK\limits_{n=1}^{\infty}\left(\frac{4n^2}{1} \low{+} \frac{4n^2}{\varphi(s,r) - r^2}\right)} = \dfrac{1}{r^2} \int_0^{\infty} \dfrac{xe^{-x\frac{\varphi(s,r)}{r}}}{\cosh x}dx. \eqno(26)$$
\end{theorem}

{\bf Proof}.
Examine equality (4) with the substitution $s:= \dfrac{\varphi(s,r)}{r}$, where $\varphi(s,r)$ is an arbitrary real-valued function of $s$ and $r$. Then we obtain the following formula for $\varphi(s,r) > r$, $r > 0$:

$$\dfrac{1}{\frac{\varphi^2(s,r)}{r^2} -1 + \TK\limits_{n=1}^{\infty}\left(\frac{4n^2}{1} \low{+} \frac{4n^2}{\frac{\varphi^2(s,r)}{r^2}-1}\right)}
=\int_0^{\infty} \dfrac{xe^{-x\frac{\varphi(s,r)}{r}}}{\cosh x}dx$$

$$\dfrac{1}{\varphi^2(s,r) - r^2 + r^2\TK\limits_{n=1}^{\infty}\left(\frac{4n^2}{1} \low{+} \frac{4n^2}{\frac{\varphi^2(s,r)}{r^2}-1}\right)}
=\dfrac{1}{r^2}\int_0^{\infty} \dfrac{xe^{-x\frac{\varphi(s,r)}{r}}}{\cosh x}dx. $$

Apply the equivalence transform with the parameters $r_0 = 1$, $r_n = r^2$, $n = 1, \ 2, \ \ldots$ to the continued fraction on the left-hand side. This results the formula:

$$\dfrac{1}{\varphi^2(s,r) - r^2 + \TK\limits_{n=1}^{\infty}\left(\frac{4n^2r^4}{r^2} \low{+} \frac{4n^2r^4}{\varphi^2(s,r)-r^2}\right)}
=\dfrac{1}{r^2}\int_0^{\infty} \dfrac{xe^{-x\frac{\varphi(s,r)}{r}}}{\cosh x}dx. $$

Using simple calculations:
$$\dfrac{1}{\varphi^2(s,r) - r^2 + \TK\limits_{n=1}^{\infty}\left(\frac{4n^2r^2}{1} \low{+} \frac{4n^2r^2}{\varphi^2(s,r)-r^2}\right)}
=\dfrac{1}{r^2}\int_0^{\infty} \dfrac{xe^{-x\frac{\varphi(s,r)}{r}}}{\cosh x}dx. $$\hfill $\square$

Let us prove the following lemma, which describes
the derivative of the continued fraction $$f(s,r) = \dfrac{1}{\varphi(s,r) + \K\limits_{n=1}^\infty \left(\frac{n^2r^2}{\varphi(s,r)}\right)}.$$

\begin{lemma} Let $\varphi(s,r) = s + \psi(r)$, where $\psi(r)$ is an arbitrary real-valued function of $r$. Then for $r>0$, $s > r - \psi(r)$
$$\frac{\partial}{\partial s}f(s,r) = - \dfrac{1}{\varphi^2(s,r) - r^2 + \TK\limits_{n=1}^{\infty}\left(\frac{4n^2r^2}{1} \low{+} \frac{4n^2r^2}{\varphi^2(s,r)-r^2}\right)}, $$ where $$f(s,r) = \dfrac{1}{\varphi(s,r) + \K\limits_{n=1}^\infty \left(\frac{n^2r^2}{\varphi(s,r)}\right)}.$$
\end{lemma}

{\bf Proof}.
Using Corollary 2, we obtain the equality $$f(s,r) = \dfrac{1}{\varphi(s,r) + \K\limits_{n=1}^\infty \left(\frac{n^2r^2}{\varphi(s,r)}\right)} = \frac{1}{r}\int_0^{+\infty} \dfrac{e^{-x\frac{\varphi(r,s)}{r}}dx}{\cosh x}.$$ Differentiating this equality by $s$ and changing the sign, we obtain: $$- \frac{\partial}{\partial s}f(s,r) = \frac{1}{r^2}\int_0^{+\infty} \dfrac{xe^{-\frac{\varphi(r,s)}{r}x}dx}{\cosh x},$$ which exactly coincide with the right-hand side of (26).\hfill $\square$

\begin{corollary} For $r> 0$, $s > \max(1, \ 2r-1)$
$$\frac{\partial}{\partial s}f_1(s,r) = - \dfrac{1}{2(1-2r+s)(1+s) + 2\TK\limits_{n=1}^{\infty}\left(\frac{4n^2r^2}{1} \low{+} \frac{4n^2r^2}{(1-2r+s)(1+s)}\right)}, $$ where $$f_1(s,r) = \dfrac{1}{2-2r+2s + 2\K\limits_{n=1}^\infty \left(\frac{n^2r^2}{1-r+s}\right)}.$$
$$\frac{\partial}{\partial s}f_2(s,r) = - \dfrac{1}{2(2r-1+s)(s-1) + 2\TK\limits_{n=1}^{\infty}\left(\frac{4n^2r^2}{1} \low{+} \frac{4n^2r^2}{(2r-1+s)(s-1)}\right)}, $$ where $$f_2(s,r) = \dfrac{1}{2r-2+2s + 2\K\limits_{n=1}^\infty \left(\frac{n^2r^2}{r-1+s}\right)}.$$
\end{corollary}

{\bf Example.} Put $\varphi(s,r) = s + \sin r$, $r = \dfrac{\pi}{2}$. Then for $s > \dfrac{\pi}{2} - 1$ we have
$$f'(s) = - \dfrac{1}{(s+1)^2 - \frac{\pi^2}{4} + \TK\limits_{n=1}^{\infty}\left(\frac{n^2\pi^2}{1} \low{+} \frac{n^2\pi^2}{(s+1)^2-\frac{\pi^2}{4} }\right)}, $$ where $$f(s) = \dfrac{2}{2s+2 + \K\limits_{n=1}^\infty \left(\frac{n^2\pi^2}{2(s+1)}\right)}.$$

\begin{theorem} For $s> \max(1, \ 2r-1)$, $r> \dfrac{1}{2}$
$$  \dfrac{\partial^2}{\partial s^2}(\ln y)(s,r) =  - \dfrac{1}{2r^2}\int_0^{\infty} \dfrac{x(e^{-\frac{1-r+s}{r}x}+ e^{-\frac{r-1+s}{r}x})}{\cosh x}dx = - h_1(s,r) - h_2(s,r), $$
where
$$h_1(s,r) = \dfrac{1}{2(1-2r+s)(1+s) + 2\TK\limits_{n=1}^{\infty}\left(\frac{4n^2r^2}{1} \low{+} \frac{4n^2r^2}{(1-2r+s)(1+s)}\right)},$$
$$ h_2(s,r) = \dfrac{1}{2(2r-1+s)(s-1) + 2\TK\limits_{n=1}^{\infty}\left(\frac{4n^2r^2}{1} \low{+} \frac{4n^2r^2}{(2r-1+s)(s-1)}\right)}.$$
\end{theorem}

{\bf Proof}.
The proof comes out from Equality 10 and Corollary 5.\hfill $\square$



{\bf Acknowledgements}. The author thanks Prof. S. Khrushchev for helpful suggestions and valuable comments.

\end{document}